\documentclass[12pt]{amsart}
\def\checkbox{\leavevmode\vbox to 9pt{\hrule \vss
	\hbox to 9pt{\vrule height 9pt \hfil\vrule height 9pt}\vss
	\hrule}\ }

\newcommand{\R}{\mathbb R}
\newcommand{\Z}{\mathbb Z}
\newcommand{\C}{\mathbb C}
\newcommand{\Tr}{{\rm Tr}}

\renewcommand{\epsilon}{\varepsilon}
\renewcommand{\phi}{\varphi}

\newtheorem{Lemma}{Lemma}[section]
\newtheorem{Theorem}{Theorem}[section]
\newtheorem{Proposition}{Proposition}[section]

\newtheorem{Corollary}{Corollary}[section]
\newtheorem{Definition}{Definition}[section]

\newtheorem{Remark}{Remark}[section]

\begin{document}

\address{Department of Mathematics, Pennsylvania State University, University Park, PA $16802$, USA
}
\author{Alexandr Borisov}
\title{Convolution Structures and Arithmetic Cohomology}
\email{borisov@math.psu.edu}
\date{\today}
\maketitle

\centerline{email: borisov@math.psu.edu}
%%%%%%%%%%%%%%%%%%%
%\begin{abstract}
\thispagestyle{empty}

%{\bf {\Large Convolution Structures and Arithmetic Cohomology }}

%\vskip .7cm

%{\sc ALEXANDR BORISOV}

%\vskip .5cm

%Department of Mathematics, Penn State University, University Park, PA $16802$, USA

%e-mail: borisov@math.psu.edu

%\vskip 1.3cm
\centerline{Revised version}

\setcounter{section}{-1}

\tableofcontents

%%%%%%%%%%%%%%%%%%%%%%%%%%%%%%%%%%%%%
\section{Introduction}

The analogy between fields of meromorphic functions on algebraic curves and global number fields has a very long history. It has been a rich source of inspiration, especially for number theory. One of the most classical objects of the theory of algebraic curves is the group of divisors. The corresponding object in the number field case is the group of Arakelov divisors (see, e.g. section 5 of this paper). In the geometric case one can associate to any divisor $D$ a linear vector space $H^0(D),$ the space of its global sections. Its dimension $h^0(D)$ plays an important role in many considerations. The classical Riemann-Roch theorem tells that for every divisor $D$ on a complete smooth curve $C$
$$h^0(D)-h^0(K-D)=\deg D +1-g,$$
where $K$ is a divisor class of differential forms and $g$ is the genus of $C$. A principal corollary of the Riemann-Roch theorem is that when $\deg D$ is big $h^0(D)=\deg D +1-g$. The arithmetic analog of this is relatively easy to obtain, see e.g. Szpiro's introduction article on Arakelov geometry (\cite{Szpiro}). To get the exact arithmetic analog of the classical Riemann-Roch formula is considerably harder. This was probably first obtained by John Tate in his thesis (cf. \cite{Tate}). 

My attention to this topic was brought by the beautiful 1998 preprint of Gerard van der Geer and Ren\'e Schoof (cf. \cite{GS}). Elaborating on the ideas of Tate, they went further to define an analog of the theta divisor, and put together a lot of information to support their choice of $h^0(D)$.

There is however one important analogy that non one was able to construct before. Namely, in the function field case a more modern version of the Riemann-Roch theorem is available, due to Serre. Roughly speaking it is the following. Besides $H^0(D)$ one can associate to $D$ another space, $H^1(D)$ (the first cohomology group) such that

1) $H^1(D)$ is dual to $H^0(K-D)$ (Serre's duality)

2) $h^0(D)-h^1(D)=deg D +1-g$ (here $h^1(D)=dim(H^1(D))$

From the first statement it follows that $h^1(D)=h^0(K-D)$ and then the classical Riemann-Roch follows from the second statement. The main goal of this paper is to construct an arithmetic analog of this  Serre's Riemann-Roch theorem.

In order to do this, we have to abandon the category of abelian groups, and use some group-like objects, the convolution of measures structures. In order to define them, one needs some notions and results from abstract harmonic analysis. The main features of our theory are the following.

1) $H^1$ is defined by a procedure very similar to \^Cech cohomology.

2) We get separately Serre's duality and Riemann-Roch formula without duality.

3) We get the  duality of $H^0(L)$ and $H^1(K-L)$ as Pontryagin duality of convolution structures.

4) The Riemann-Roch formula of Tate - van der Geer-Schoof follows automatically from our construction by an appropriate dimension function.

The paper is organized as follows. In section 1 we recall some necessary definitions and results from harmonic analysis. In section 2 we define our basic objects (ghost-spaces) and their dimensions.  In section 3 we introduce some short exact sequences of ghost-spaces. In section 4 we develop the duality theory of ghost-spaces. In section 5 we apply the theory to arithmetic and obtain our main results. In section 6 we discuss possible directions in which the theory can grow.

{\bf Acknowledgments.} The author thanks Adrian Ocneanu, Yuri Zarhin, and Nik Weaver for their interest and stimulating discussions. The author is especially thankful to Michael Voit for the expert's advises on convolution structures. The author also thanks Jeff Lagarias for the numerous helpful comments on the preliminary draft of the paper.

%%%%%%%%%%%%%%%%%%%%%%%%%%%%%%%%%%%%%%%%%%%%%%%
\section{Some results from harmonic analysis}

This section is intended primarily for arithmetic geometers and other readers with little knowledge of harmonic analysis. All of the material presented here is contained in basic harmonic analysis textbooks. We claim no originality whatsoever.

Most proofs are omitted. All the missing proofs can be found, e.g. in the book of Folland \cite{Folland}. See also Berg and Forst \cite{BergForst} for a more in-depth treatment of positive-definite functions. We are only interested in the commutative case even though many of the results are true for arbitrary locally compact groups. We start with some basic definitions.

\begin{Definition} A locally compact abelian group $G$ is a topological group which is abelian and locally compact. We will use addition notation for the operation, unless $G\subseteq \C^*$.
\end{Definition}

\begin{Definition} A character $\chi$ on a locally compact abelian group $G$ is a group homeomorphism $\chi : G \rightarrow \C^*$.
\end{Definition}

\begin{Definition} A character $\chi$ is called symmetric if $\chi(-x)=\overline{\chi (x)}.$ This is equivalent to the image of $\chi$ being a subset of the unit circle.
\end{Definition}

\begin{Definition} All characters of $G$ form an abelian group, under the pointwise multiplication. Its subgroup consisting of symmetric characters is called the Pontryagin dual group of $G$. It can be endowed with a natural topology of pointwise convergence. It is denoted by $\widehat{G}$.
\end{Definition}

The following theorem is well-known in abstract harmonic analysis.

\begin{Theorem} (Pontryagin Duality) The group $\widehat{G}$ is locally compact. Its Pontryagin dual $\widehat {\widehat{G}}$ is naturally isomorphic to the group $G$. 
\end{Theorem}

The isomorphism in the above theorem is the following. Every $x\in G$ gives a function on $\widehat{G}$ by sending $\chi$ to $\chi(x)$. This defines a map from $G$ to $\widehat {\widehat{G}}$. It is easy to show that this map is a homeomorphism. The above theorem states that it is an isomorphism.

Some of the main instruments in the proof of the Pontryagin Duality Theorem, which we will also use a lot, are the notions of the Haar measure, the Fourier transform and the inverse Fourier transform.
Here are some basic definitions and results.

\begin{Theorem} Suppose $G$ is a locally compact group (not necessarily abelian). Then there exists a non-zero left-invariant $\sigma -$additive Borel measure on it. This measure is unique up to a multiplicative constant. It is called a left Haar measure. When $G$ is abelian (or, in general, if this measure is right-invariant as well) it is called simply a Haar measure of the group $G.$ 
\end{Theorem}

\begin{Definition} Suppose $m_G$ is a Haar measure on a locally compact abelian group $G.$ Suppose $f$ is a function on $G$. Then the Fourier transform of $f$, relative to the measure $m_G$ is the function on $\widehat{G}$ defined as follows.
$$\hat{f}(\chi )= \int \limits_{x\in G} f(x) \overline{\chi (x)} dm_G(x) $$
\end{Definition}

The Fourier transform is defined for all functions on $G$ that are $L^1$ with respect to a Haar measure.

\begin{Definition} Suppose $\mu$ is some complex-valued measure on $\widehat{G}.$ Then its inverse Fourier transform $\check{\mu}$ is a function on $G$ defined as follows.
$$\check{\mu}(x)= \int \limits_{\chi \in \widehat{G}} \chi(x) d\mu (\chi )$$
\end{Definition}

The inverse Fourier transform is defined for all bounded measures $\mu$ on $\widehat{G}.$ Unlike a Fourier transform it does not involve a choice of a Haar measure.

We are now going to discuss the notions of positive-definite functions and measures on locally compact abelian groups.

\begin{Definition} A complex-valued function $f$ on a group $G$ is called positive-definite iff for all $x_1,x_2,\dots , x_n \in G$ the matrix $f(x_j-x_i)$ is hermitian nonnegative-definite.
\end{Definition}

Note that a positive-definite function need not be continuous. Examples of positive-definite functions include the characteristic functions of subgroups of $G$. Another important example is a function  $e^{-Q(x)}$ on $\R ^n$ where $Q(x)$ is a positive-definite quadratic form. For more examples see \cite{BergForst}, chapter 1, section 5.

Positive-definite functions have many interesting properties, some of which will be discussed later. One of the most important results about them, which can be viewed as a step toward the Pontryagin Duality Theorem, is the following theorem of Bochner.

\begin{Theorem} Suppose $f$ is a continuous positive-definite function on a locally compact abelian group $G.$ Then there exists a unique measure $\mu$ on $\widehat{G}$ such that $f=\check{\mu}.$ This measure is real-valued and non-negative.
\end{Theorem}

For any topological space $G$ one can multiply functions on $G$. If  $G$ is given a structure of a locally compact group, then we get additionally an operation of convolution of measures. This operation can be described as follows. Given two bounded measures $\mu $ and $\nu$ on $G$, one can take their Cartesian product, which is a measure on $G\times G$. Then their convolution $\mu *\nu $ is the pushforward of that product with respect to the addition map $G\times G \rightarrow G.$ One can also define the convolution of functions by the following integral.
$$(f*g)(x)=\int \limits_{y\in G} f(y)g(x-y) dm(y)$$
Here $m$ is a Haar measure on $G$. This agrees with the operation of convolution of measures in a natural way:
$$(f\cdot m)*(g\cdot m)=(f*g)\cdot m$$
Using the above convolution we can make the following definition.

\begin{Definition} A measure $\mu$ is called positive-definite iff for any continuous function $f$ on $G$ with compact support
$$\int \limits_G (f * \bar{f}) d\mu \ge 0$$
\end{Definition}

The point measure at $0$ and the Haar measures are always positive-definite. If $f$ is a positive-definite function, which is $L^1$ with respect to a Haar measure $m$, then $f\cdot m$ is a positive-definite measure. Also any measure whose inverse Fourier transform is real-valued and nonnegative, is positive-definite.

The convolution of measures makes the space of bounded measures on a locally compact group an algebra. One can recover the group structure on the set $G$ from the operation of convolution of measures. Namely, the convolution of two point measures is a point measure of the sum. Pontryagin duality essentially switches the algebra of functions and the algebra of measures, via the Fourier transform. This suggests that one can generalize the notion of a locally compact group to a arbitrary sets $X$ together with some algebra structure on some subspace $M$ of the space of measures on $X.$ Depending on algebraic and analytic restrictions on this convolution algebra, many different versions of this were proposed. For some of those classes of structures the Pontryagin duality theorem holds. For a quite general framework and a good survey, see \cite{Vain}. Harmonic analysts are mostly interested in the noncommutative situation. We only need commutative convolution structures.

\begin{Definition} Suppose $G$ is a topological space. The {\bf weak topology} on the space of measures on $G$ is the weakest topology such that for every continuous function with compact support on $G$ the corresponding linear operator is continuous.
\end{Definition}

\begin{Definition} A commutative convolution of measures structure $*$ on a space $G$ is called weakly separately continuous if and only if  for any measure $\mu \in M$ the linear operator from $M$ to $M$ sending $\nu$ to $\mu * \nu$ is weakly continuous in $\nu$.
(If the convolution is not necessarily commutative then one also needs $\nu * \mu$ to be weakly continuous in $\nu$)
\end{Definition}

The linear combinations of point measures are dense in the set of all measures with respect to the weak topology. Because of this any weakly separately continuous convolution of measures structure is uniquely determined by what it does on the point measures (see Pym \cite{Pym} for a more detailed discussion of this). Thus we can identify the convolution structures with the map $* :G\times G \rightarrow Measures(G)$ given by
$$(x,y)\mapsto \delta_x*\delta_y$$
This will be our convention in the next section. Of course, not all maps as above extend to associative operations, so this always has to be checked.
%%%%%%%%%%%%%%%%%%%%%%%%%%%%%%%%%%%%%%%%%%%%%%%
\section{Ghost-spaces and their dimensions}

First we would like to explain our motivation.

Suppose $G$ is a locally compact abelian group. Suppose $H$ is a subgroup of $G$. Then its characteristic function $i_H$ has the following properties.
\begin{enumerate}
\item $i_H(0)=1$

\item $i_H$ is even (i.e. $i_H(-x)=i_H(x)$ for all $x\in G$)

\item $i_H$ is positive-definite

\item $(i_H)^2=i_H$
\end{enumerate}

One can check that any function satisfying the four conditions above is a characteristic function of some subgroup of $G$.

In their paper \cite{GS} van der Geer and Schoof defined $h^0(D)$ as a logarithm of the sum of $e^{-Q(x)}$ over some lattice, where $Q(x)$ is a quadratic form on that lattice. They did not formally define $H^0(D)$ but they essentially viewed it as a ``subgroup" of the lattice defined by a ``characteristic function" $e^{-Q(x)}.$ Note that this function satisfies all of the conditions above except the last one. So we want to generalize the notion of a subgroup (and therefore a group) by abandoning this last condition. In fact, the first three conditions  already have some interesting implications. For example, they imply that the corresponding function never takes values bigger than $1$. See Theorem 2.1 below for the proof.

Because we want our sub-object to be supported on the whole group $G$ we will also assume that its ``characteristic function" is strictly positive. We will also assume that it is continuous, because this is what we have in applications. We would like to consider this sub-object as a space with a convolution of measures structure. We will call such a sub-object a ghost-space (because its elements ``only exist with some probability"). More precisely, we will call it a ghost-space of the first kind, in order to distinguish it from another kind of convolution structures, that will be defined later.
One can also unify the two kinds of ghost-spaces. This more general kind of objects (see Example in section 6) will probably be needed in order to extend the theory to higher dimensions. In this paper we will only develop the theory of ghost-spaces to the minimal extent necessary for the applications to the number fields.

\begin{Lemma} Suppose $G$ is a locally compact abelian group. Suppose $u:G\rightarrow R^+$ is a positive, positive-definite continuous function on it such that $u(0)=1.$ Consider the convolution of measures $*$ on $G$ such that
$$\delta_x *\delta_y =\frac{u(x)u(y)}{u(x+y)}\delta_{x+y}$$
Then this convolution is commutative and associative.
\end{Lemma}

{\bf Proof.} Consider the space of all measures $\mu $ with the property that $u\cdot \mu$ is bounded. We can make it a convolution algebra by setting
$$\mu_1 * \mu_2 =\frac{(u\mu_1)\odot (u\mu_2)}{u},$$
where $\odot$ is the standard convolution of measures on $G.$ This convolution $*$ extends the convolution $\delta_x *\delta_y.$ It is obviously commutative, and associative. It is also weakly separately continuous, where the weak topology is defined using the continuous functions with compact support.
\hfill \checkbox

\begin{Definition}  We will call the pair $(G, *),$ where $*$ is a non-standard convolution on $G$ as above, the ghost-space of the first kind. We will denote it by $G_u$.
\end{Definition}

\begin{Remark}
The function $u$ can be recovered from the convolution $*$ as its only real-valued positive symmetric quasi-character. See section 4 for the details.  See also Voit (\cite{Voit}) for a related more general theory.
\end{Remark}

\begin{Remark} One can see that the convolution algebras of $G_u$ and $G$ are isomorphic. However the pair of algebras (functions and measures) on $G_u$ together with the (function, measure) pairing is different.
\end{Remark}

\begin{Definition} We define the dimension of $G_u$ which depends on the choice of a Haar measure $m$ on $G$ as follows.
$$\dim _m G_u = \log \int \limits_G u(x) dm(x) $$
When $G$ is discrete, it has a distinguished Haar measure, the counting measure $m_c$. In this case we will say that the absolute dimension of $G_u$ 
$$\dim G_u =\dim_{m_c} G_u.$$
(When $G$ is not discrete, we will think of $G_u$ as having infinite absolute dimension.)
\end{Definition}

{\bf Examples.}

1) Suppose $G$ is a locally compact abelian group. Then $G_1$ is just $G$ itself with the standard convolution of measures. We will therefore identify $G_1$ with $G.$

2) Suppose $G=\Z ^n$ and $Q$ is a positive-definite quadratic form on it. Then one can check that $u(x)=e^{-Q(x,x)}$ is positive-definite (see, e.g. \cite{BergForst}, Proposition 7.19). So one can define the ghost space $G_u.$ Its dimension, in the above sense, is equal to  $\log \sum \limits_{x\in \Z} e^{-Q(x,x)}$. This is exactly the kind of formula that van der Geer and Schoof used to define $h^0(D),$ and $u(x)$ is their effectivity function. So the finite-dimensional ghost-space of the first kind $G_u$ is going to be, in our interpretation, $H^0(D).$

The following easy theorem supports our interpretation of $u$ as a functions that measures the ``probability with which elements of $G_u$ exist".

\begin{Theorem} Suppose $G_u$ is a ghost-space of the first kind. Then for all $x\in G$ $u(x)\le 1$. Also, those $x$ that $u(x)=1$ form a closed subgroup $H$ of $G.$  Moreover, $u(x)$ comes from a function on $G/H.$
\end{Theorem}

{\bf Proof.} The first claim is contained in Folland (see \cite{Folland}, cor. 3.32.) To prove the second and third claims we note that by \cite{Folland}, prop. 3.35 the following matrix is positive-definite.
$$ \left [\begin {array}{lcr} \mbox {1}&\mbox {u(x)}&\mbox {u(x+y)}
\\\noalign{\medskip}\mbox {u(x)}&\mbox {1}&\mbox {u(y)}
\\\noalign{\medskip}\mbox {u(x+y)}&\mbox {u(y)}&\mbox {1}\end {array}
\right ]
 $$
If $u(x)=1$, it implies that $(u(x+y)-u(y))^2\le 0,$ so $u(x+y)=u(y).$ This implies the theorem. \hfill \checkbox

\vskip 0.2cm

Now we define the ghost-spaces of the second kind. While the ghost-spaces of the first kind are intuitively the abelian groups with ``partially existent" elements, the ghost-spaces of the second kind have different nature. Their elements exist with the probability $1$, but their position on $G$ is not fixed. They could be thought of as ``clouds" on $G$. As a result the ``addition" of two such elements is probabilistic. More precisely, the addition will have a translation-invariant error probability. This kind of objects appears in particular when one tries to take a quotient of $G$ over its sub-object $G_u$. We refer to the next section for a more detailed explanation. Now we just give a formal definition.

\begin{Definition} Suppose $G$ is a locally compact abelian group. Suppose $\mu$ is a positive-definite, even, positive probability measure  on $G.$ We will call the pair $(G,*)$ with the convolution of measures $*$ from the next lemma the ghost-space of the second kind. It will be denoted by $G^{\mu}.$
\end{Definition}

\begin{Lemma}
Suppose $G$ and $\mu$ are as above. Consider the convolution of measures $*$ on $G$ such that
$$\delta_x *\delta_y = T_{x+y}\mu ,$$
where  $T_{x+y}$ is the usual shift by $(x+y)$. Then this convolution is commutative and associative.
\end{Lemma}

{\bf Proof.} We will show that $*$ extends to the space of bounded measures. We will use for that the canonical continuation formula of Pym (cf. \cite{Pym}). For any two bounded  measures $\nu_1$ and $\nu_2,$ and a continuous function with compact support $f$ on $G,$ the following formula makes sense.
$$(\nu_1 * \nu_2) (f) = \int \int (T_{x+y}\mu )(f) d\nu_1(x) d\nu_2(y)$$
One can use it to define the  measure $\nu_1 * \nu_2$. This obviously generalizes the convolution $*$ from the statement of the lemma. One can easily check that $\nu_1 * \nu_2$ is bounded. Moreover, the convolution of two probability measures is a probability measure, and the convolution is weakly separately continuous. We now need to check that it is associative.
If $\nu_1,$ $\nu_2,$ $\nu_3$ are bounded  measures and $f$ is a continuous function with compact support on $G$ then one can check the following.
$$\Big( (\nu_1 *\nu_2) * \nu_3 \Big) (f) = \int \int \int \Big( T_{x+y+z}(\mu \odot \mu) \Big) (f) d\nu_1(x) d\nu_2(y) d\nu_3(z)  ,$$
where $\odot $ is the standard convolution of measures on $G.$ The associativity follows.
\hfill \checkbox

Obviously, the measure $\mu$ is uniquely determined by $*$.
Also if $\mu = \delta_0$ then the convolution above is just a standard convolution of measures on $G$. Thus we will identify $G^{\delta_0}$ with $G$.

\begin{Definition} Suppose $G^{\mu}$ is a ghost-space of the second kind. Suppose $m$ is a Haar measure on $G$. Suppose $\mu$ is absolutely continuous with respect to $m$, i.e. $\mu = u \cdot m$ for some function $u$ on $G$. Then we define
$$ \dim ^{(m)} G^{\mu} = \log u(0)$$
(If $\mu$ is not absolutely continuous, we will think of $G^{\mu}$ as being infinite-dimensional. In this paper we only consider the finite-dimensional ghost-spaces of the second kind).

If $G$ is compact, then it has a distinguished Haar measure, the probability measure $m_{prob}$. In this case we define the absolute dimension
$$ \dim G^{\mu} = \dim ^{(m_{prob})} G^{\mu} $$

\end{Definition}

Some justification of the above definition is provided by Lemma 2.3. The real justification, however, is in Proposition 3.1 and Theorem 4.1.

One can choose to consider a group $G$ itself both as a ghost-space of the first and of the second kind. In fact, one can see immediately that this is the only case when a convolution structure can be interpreted in these two ways. In order for our notation to be consistent we need to check that the absolute dimension of $G$ is does not depend on this interpretation. If $G$ is either not discrete or not compact then it has infinite absolute dimension. So the only case we really need to consider is when $G$ is finite. The following lemma does just that.

\begin{Lemma} Suppose $G$ is a finite abelian group. Then its dimension as a ghost-space of the first or the second kind is equal to $\log |G|.$
\end{Lemma}

{\bf Proof.} We will denote by $M$ the counting measure on $G.$

1) As a ghost-space of the first kind $G=G_1.$ So $\dim G =\dim_M G= \log |G|.$
\hfill \checkbox

2) As a ghost-space of the second kind $G=G^{\delta_0}.$ If $m$ is the probability Haar measure on $G,$ then $m=\frac{1}{|G|}M.$ So $\delta_0 =h\cdot m,$ where $h(0)= |G|,$ $h(x)=0 $ for $x\ne 0.$ Therefore $\dim G =\dim G^{\delta_0}=\log h(0) =\log |G|.$ \hfill \checkbox

%%%%%%%%%%%%%%%%%%%%%%%%%%%%%%%%%%%%%%%%%%%%%%%%%%%%%%%%%%%%%%%%%
\section{Short exact sequences of ghost-spaces}
In this section we will define some short exact sequences of ghost-spaces. We will check that the dimension is additive, whenever defined. We must note that this is probably just a little piece of the more general theory which is yet to be developed.

First of all, we want to consider $G_u$ as a sub-object of $G$. Let's try to define a quotient object $G/G_u$. We want our definition to be roughly parallel to the definition of the group quotient $G/H$ where $H$ is a closed subgroup of $G$. In this latter situation the objects of the quotient space can be identified with the cosets of $H.$ If $H$ is compact, one can associate to each coset $x+H$ the probability measure $T_x(p_*(m_{prob}(H))).$ Here $m_{prob}(H)$ is the Haar probability measure on $H,$ $p:H\rightarrow G$ is the embedding, and $T_x$ is translation by $x$. Then the convolution of the measures corresponding to $x+H$ and $y+H$ is the measure that corresponds to $x+y+H$.

Suppose now that $G_u$ is a ghost space, with $\dim _m G_u <\infty$ for some Haar measure $m$ on $G$. The natural analogs of the measures above are $T_x (\mu )$ where $\mu$ is the probability measure proportional to $u\cdot m.$ When one convolves two such measures using the standard convolution on $G$ (that corresponds to the addition in the ambient group) one gets the following
$$T_x(\mu) * T_y(\mu) = T_{x+y}(\mu *\mu)= \int \limits_{z\in G}(T_{x+y+z}\mu)d\mu (z)$$
Thus the quotient $G/G_u$ can be viewed as a ghost-space of the second kind $G^{\mu}.$ Here is the formal definition.

\begin{Definition} Suppose $G_u$ is a ghost-space of the first kind. Then we say that $G_u$ is a subspace of $G.$  If $\dim G_u <\infty$ we also say that the quotient $G/G_u$ is the ghost-space of the second kind $G^{\mu},$ where $\mu$ is the probability measure on $G$ proportional to $u(x)\cdot m$. Here $m$ is some (any) Haar measure on $G.$
\end{Definition}

\begin{Proposition} The dimension is additive in the above short exact sequence, provided we use the same Haar measure for $G$ and $G_u$ to define it. That is, whenever defined,
$$\dim _m G= \dim _m G_u + \dim G^{\mu}.$$
\end{Proposition}

{\bf Proof.} Because of the Definition 2.3 we only need to consider the case when $G$ is compact. Since changing the Haar measure $m$ has no effect on the validity of the above identity, we can choose $m$ to be the probability measure. If $\dim_m G_u = \log A$ then $\mu = \frac{1}{A}\cdot u\cdot m.$ Therefore $\dim_m G =0,$ $\dim_m G_u =\log A,$ and  $\dim G^{\mu} = \log (\frac{u(0)}{A})=-\log A.$ The last identity is because $u(0)=1$ by the definition.
\hfill \checkbox

\vskip .2cm

Now we define another kind of short exact sequences. This time all objects are ghost-spaces of the first kind.

\begin{Definition} Suppose $G$ is a locally compact abelian group and $H$ is its closed subgroup. Suppose $u:G\rightarrow R^+$ is a positive-definite, positive, even, continuous function on $G$ such that $u(0)=1.$ Abusing notation a little bit, we will call the restriction of $u$ to $H$ also $u.$ Then we will say that $H_u$ is a subspace of $G_u$. If we can define a positive-definite continuous function of $v$ on $G/H$ as below we will also say that $(G/H)_v$ is the quotient $G_u/H_u$.
$$v(xH)=\frac{\int\limits_{y\in H}u(x+y)dm(y)}{\int\limits_{y\in H} u(y)dm(y)},$$ 
where $m$ is a Haar measure on $H.$
\end{Definition}

\begin{Remark} In fact, $v$ is probably always positive-definite, whenever it is defined and continuous. At least it is true if both $\dim G$ and $\dim H$ are finite, as the following proposition shows.
\end{Remark}

\begin{Proposition} Suppose $u$ and $v$ are continuous functions defined as in Definition 3.2. Suppose that $\int \limits_G u(x) dm_G(x)$ and $\int \limits_H u(x) dm_H(x) $ are both finite. Then $v$ is positive-definite.
\end{Proposition}

{\bf Proof.} Since $v \in L^1(G/H)$ it is enough to show (cf. \cite{Folland}, 4.17) that
$$\int \limits_{G/H} \chi(y) v(y) dm_{G/H}(y) \ge 0$$
for any character $\chi $ on $G/H$. By the definition of $v$ it is equivalent to saying that 
$$\int \limits_{G} \chi(x) v(x) dm_{G}(x) \ge 0$$
for all characters $\chi $ on $G$ that come from $G/H.$ This now follows from $u$ being continuous and positive-definite (cf. \cite{Folland}, 4.23).
\hfill \checkbox

\begin{Remark} The dimension is obviously additive in the above short exact sequence if one chooses the measure on the quotient space as the quotient of measures on $G$ and $H.$
\end{Remark}

\begin{Remark} Pretty obviously, $G_1/H_1=(G/H)_1$ whenever defined (i.e. when $H$ is compact). So our definition really is compatible with the usual group quotients. 
\end{Remark}

\begin{Remark} One can also define similarly some short exact sequences of the ghost-spaces of second kind. They will be dual to the above short exact sequences in the sense of the next section.
\end{Remark}

%%%%%%%%%%%%%%%%%%%%%%%%%%%%%%%%%%%%%%%%%%%%%%%
\section{Duality theory of ghost-spaces}

Here we develop the duality theory of ghost-spaces. Basically, the dual of $G_u$ is $\widehat{G}^{\hat{u}},$ where  $\widehat{G}$ is the Pontryagin dual of $G$ and $\hat{u}$ is the Fourier transform of $u.$ To be precise,  $\hat{u}$ is such measure that
$$u(x)=\int\limits_{y\in \widehat{G}} y(x) d\hat{u}(x).$$
The existence of such measure is the Bochner theorem on $G$ (cf., e.g. Folland \cite{Folland}, prop. 4.18). We could have taken this as a definition, of course. But we already had a lot of ad hoc definitions in the previous two sections. So we claim that this duality really is the Pontryagin duality of convolution structures.

We should mention here that a lot of work has been done by researchers in harmonic analysis to extend Pontryagin duality of locally compact abelian groups to the more general convolution structures. We should mention here for reference the survey of Vainerman \cite{Vain}. It looks like the particular case we need is new. But it is very similar algebraically to the more general case  of commutative signed hypergroups, as introduced by Margit R\"osler ( \cite{Roesler1}, \cite{Roesler2} ).  To be precise, for any $G_u$ one can define an involution by sending $x$ to $-x,$ and a measure $\omega =\frac{m}{u^2},$ where $m$ is some Haar measure on $G.$ Then the triple $(G,\omega, *)$ satisfies the algebraic part of the axioms of a commutative signed hypergroup.

So we will construct the dual of $G_u$ following the construction of R\"osler. We are only interested in the algebraic part of the construction, and our convolutions are given by explicit formulas. So we will basically ignore the analytic part of the theory.

First, let us consider all quasi-characters on $G.$ These are the functions $\phi : G\rightarrow \C$ with the following property.
$$\phi(x) \cdot \phi(y) =\int \limits_G \phi(\lambda )(\delta_x *\delta_y) (\lambda )$$
In our case this means that
$$\phi(x) \cdot \phi(y) =\phi (x+y) \frac{u(x)u(y)}{u(x+y)}$$
So $\frac{\phi(x)}{u(x)}$ is a multiplicative function on $G.$ This implies that $\phi(x)=\chi(x)u(x)$ for some multiplicative function $\chi : G\rightarrow \C$.

Now we should consider only the symmetric quasi-characters, i.e. those $\phi$ that $\phi(-x)=\overline{\phi(x)}.$ One can see from the above description of quasi-characters that these are $\phi_{\chi}(x) =\chi(x)u(x) $ for some $\chi : G\rightarrow S^1,$ i.e. for $\chi \in \widehat{G}.$

So we established the natural set-wise isomorphism of $\widehat{(G_u)}$ and $\widehat{G}.$ We can therefore transfer the group structure of $\widehat{G}$ onto $\widehat{(G_u)}.$ What we really need to do though is to figure out the convolution structure on $\widehat{(G_u)}$. First we can define the Fourier transform and the inverse Fourier transform as in R\"osler \cite{Roesler1}.

Since $\phi_{\chi}(x) =\chi(x)u(x) $, for all $x\in G$, we have that
$$\check{\delta_{\chi}}(x)=\chi(x)u(x),$$
where $\delta_{\chi}$ is a point measure at $\phi_{\chi}$.

The convolution of measures in  $\widehat{(G_u)}$ should correspond via the inverse Fourier transform to the multiplication of functions on $G_u,$ i.e. to the usual multiplication of functions on $G$. The only thing we really need to prove is the  following proposition.

\begin{Proposition} Suppose $\chi_1, \chi_2 \in \widehat{G},$ $x\in G$. Then
$$(\chi_1(x)u(x))\cdot (\chi_2(x)u(x))= \int \limits_{\chi \in \widehat{G}} \chi(x)u(x) d(T_{\chi_1+\chi_2}\hat{u})(\chi)$$
\end{Proposition}

{\bf Proof.} The above equality is equivalent to the following.
$$u(x)=\int \limits_{\chi \in \widehat{G}} \frac{\chi (x)}{\chi_1(x)\chi_2(x)} d(T_{\chi_1+\chi_2}\hat{u})(\chi)$$
The right hand side can be rewritten as 
$$\int \limits_{\chi \in \widehat{G}} (\chi -\chi_1-\chi_2)(x) d(T_{\chi_1+\chi_2}\hat{u})(\chi) $$
Using the substitution $\lambda =\chi -\chi_1-\chi_2,$ it is equal to

$$\int \limits_{\lambda \in \widehat{G}} \lambda (x) d\hat{u}(\lambda)$$
Then the desired equality is just the definition of $\hat{u}.$  \hfill \checkbox

One can also check that the natural involution of quasi-characters $\phi \mapsto \bar{\phi} $ corresponds to $\chi \mapsto -\chi.$ To complete the picture we need to show that $\widehat{\widehat{(G_u)}}$ is naturally isomorphic to $G_u.$ This means that all the symmetric quasi-characters of the convolution structure $\widehat{G}^{\hat{u}}$ are of the form $\chi(x) u(x)$ for some $x \in G.$
The following proposition does just that.

\begin{Proposition} Suppose $f: G\rightarrow \C$ is a symmetric quasi-character on $\widehat{G}^{\hat{u}}.$ Then $f(x)=\chi(x) u(x)$ for some $x \in G.$
\end{Proposition}

{\bf Proof.} Being a quasi-character here means that for all $\chi_1, \chi_2 \in \widehat{G}$
$$f(\chi_1)\cdot f(\chi_2) = T_{\chi_1+\chi_2} \hat{u} (f).$$
Therefore 
$$f(\chi_1)\cdot f(\chi_2) = f(0) \cdot f(\chi_1+\chi_2).$$
This implies that $f(\chi)=v(\chi) \cdot f(0),$ where $v$ is a character on $\widehat{G}.$

Also, since $f$ is symmetric, $f(0)=\bar{f(0)},$ so $f(0) \in \R .$ As a result, the condition $f(-\chi) = \bar{f(\chi )}$ implies that $v(-\chi) = \bar{v(\chi)}$ so $v$ takes values in the unit circle $S^1.$ By the Pontryagin duality theorem, $v(\chi) =\chi(x)$ for some $x \in G .$

Finally, $f(0) \cdot f(0) = \hat{u} (v\cdot f(0).$ So $f(0) =\hat{u}(v).$ By the definition of $\hat{u},$ $f(0)=u(x),$ the proposition is proven. \hfill \checkbox

\begin{Remark} If we take duals in a short exact sequence of Definition 3.2 we get again a short exact sequence, going in the opposite direction. So the situation is completely parallel to the case of usual locally compact abelian groups.
\end{Remark}

Now let's discuss what happens with the dimension when the dual is taken. First of all, $\dim \widehat{G}^{\hat{u}}$ only makes sense if $\widehat{G}$ is compact, and $\hat{u}$ is absolutely continuous with respect to a Haar measure. This means that $G$ is discrete. Then we have the following theorem.

\begin{Theorem} Suppose $G$ is discrete, $G_u$ is a finite-dimensional ghost-space of the first kind. Then
$$\dim G_u=\dim \widehat{G_u}$$
\end{Theorem}

{\bf Proof.} Consider the counting measure $m$ on $G.$ Its dual measure $\hat{m}$ is a probability Haar measure on $\widehat{G}$ (cf., e.g. Folland \cite{Folland}, Prop. 4.24). Then $\hat{u} = f(\chi)\cdot \hat{m}$ where $f$ is the Fourier transform of $u$ relative to the above measures (cf. Folland, \cite{Folland}, prop. 4.21). By definition, 
$$\dim \widehat{G}^{\hat{u}} =\log f(0) = \dim G_u$$
 \hfill \checkbox

\begin{Remark}
Even though it might be possible to extend the definition of the dimension of the ghost-spaces of the second kind, the above theorem is not likely to have any generalizations. The following example highlights the major obstacle.
\end{Remark}

{\bf Example.} Suppose $u=e^{-\pi x^2}$ is a function on $\R ,$ and $m$ is the standard measure on $\R .$ Then $\R_u$ is the ghost-space of the first kind and $\R^{u m}$ is the ghost-space of the second kind. We have the following short exact sequence of ghost-spaces.
$$0 \rightarrow \R_u \rightarrow \R \rightarrow \R^{u m} \rightarrow 0$$
We have that $\dim \R =\infty $. For any measure $M$ $\dim _M \R_u$ is finite (equal to zero if $M=m$). By the nature of dimension, we expect that $\dim \R^{u m} =\infty .$ On the other hand, one can check that $\widehat{\R_u} = \R^{u m}.$ So we have a duality between a finite-dimensional ghost-space $\R_u$ and an infinite-dimensional ghost-space $\R^{u m}.$

%%%%%%%%%%%%%%%%%%%%%%%%%%%%%%%%%%%%%%%%%%%%%%%%%%%
\section{Arithmetic cohomology via ghost-spaces}

First of all, let us fix the  same notations as in  \cite{GS}, section 3. For the convenience of a reader we reproduce most of them below.

Our main object is  an ``arithmetic curve", i.e. a number field $F.$ An Arakelov divisor $D$ on it is a formal sum $\sum \limits_P x_P P+\sum \limits_{\sigma} x_{\sigma}\sigma,$ where $P$ runs over the maximal prime ideals of the ring of integers $O_F$ and $\sigma$ runs over the infinite, or archimedean places of the number field $F$. The coefficients $x_P$ are in $\Z$ while the coefficients $x_{\sigma}$ are in $\R$. The degree $\deg(D) = \sum \limits_{P} \log (N(P))x_P + \sum \limits_{\sigma} x_{\sigma} $.

An Arakelov divisor $D$ is determined by the associated fractional ideal $I=\prod P^{-x_p}$ and by $r_1+r_2$ coefficients $x_{\sigma} \in \R$. We can define a hermitian metric on $I,$ and on $I\otimes \R=F\otimes \R$ as in \cite{GS}. That is, for $z=(z_{\sigma})$
$$||(z_{\sigma})||^2_D=\sum \limits_{\sigma}|z_{\sigma}|^2\cdot ||1||^2_{\sigma},$$
where $||1||^2_{\sigma}=e^{-2x_{\sigma}}$ for real $\sigma$ and  $||1||^2_{\sigma}=2e^{-x_{\sigma}}$ for complex $\sigma$.
According to van der Geer and Schoof,
$$h^0(D)=\sum\limits_{x\in I} e^{-\pi ||x||^2_D}$$
In accordance with this, we make the following definition.

\begin{Definition} In the above notations, $H^0(D)$ is the ghost-space of the first kind $I_u,$ where $u(x)= e^{-\pi ||x||^2_D}$.
\end{Definition}

\begin{Remark}
To make the above definition valid, we need to check that $u$ positive-definite. This basically follows from the positivity of its Fourier dual, which will be calculated in Theorem 5.2 (cf., e.g. Folland \cite{Folland}). Clearly, $\dim I_u=h^0(D).$
\end{Remark}

Now we are going to define $H^1(D)$. First, let us look at how it can be done in the geometric situation. We have the curve $C$ with the map $\pi:C\rightarrow P^1.$ Probably the easiest way to calculate $H^1(D)$ in this situation is by \^Cech cohomology. For this we need to cover the curve by affine open sets. One way to do it is to choose two points on $P^1$, say $\alpha$ and $\infty$, and consider the open sets $U_0=\pi^{-1}(P^1 - \infty ) $ and $U_1 =\pi^{-1}(P^1 - \alpha ) $. Then we have the following four spaces.
$$V_{00}=H^0(D, U_0\cap U_1)$$
$$V_{10}=H^0(D, U_0) $$
$$V_{01}=H^0(D,  U_1) $$
$$V_{11}=H^0(D) $$
Here $V_{10}$ and $V_{01}$ are subspaces of $V_{00}$ and $V_{10}\cap V_{01} =V_{11}.$ By the definition of \^Cech cohomology, and since $U_0$ and $U_1$ are affine, 
$$H^1(D)=V_{00}/(V_{01}+V_{10}) =(V_{00}/V_{10})/(V_{01}/V_{11})$$
Now we try something similar in the arithmetic case. Let us choose $U_0 =\pi^{-1}(\infty)$ and $U_1=\pi^{-1}(p)$ where $p$ is some prime number. Let us denote by $J$ the localization of $I$ in $p$. Then the natural analog of $V_{11}$ above is the ghost space $I_u$ for  $u(x)= e^{-\pi ||x||^2_D}$. The analog of $V_{10}$ is $I$. The analog of $V_{00}$ is $J$. The analog of $V_{01}$ would have been $J_u$, if we managed to define ghost-spaces for the groups like $J$. Then the \^Cech cohomology of this covering should be
$$(J/I)/(J_u/I_u).$$
Now we have some problems. It looks like the different choices of $p$ should lead to different answers, unless we are willing to complete $J$ to $I\otimes \R$. So this is what we do. Please note that $I\otimes \R$ is a locally compact group, and we have no problems in defining the ghost-space $V_{01}$. We also have no problems to define other ingredients in the formula using the short exact sequences from section 3.
So this is our definition.

\begin{Definition} For an Arakelov divisor $D$ as above
$$H^1(D)= ((I\otimes \R) /I)/((I\otimes \R)_u/I_u)$$
Also, $h^1(D)=\dim H^1(D),$
as the dimension of the ghost-space of the second kind.
\end{Definition}

We will see that this definition yields a beautiful theory with such attributes of the geometric case as Serre's duality and Riemann-Roch. For this we just need to do some calculations.

\begin{Proposition} We have that 
$$(I\otimes \R)_u/I_u = ((I\otimes \R) /I)_v,$$
where for every $\overline{x} \in (I\otimes \R) /I$
$$v(\overline{x}) = \frac{\sum \limits_{y \in I} e^{-\pi ||x+y||^2_D}}{\sum \limits_{y \in I} e^{-\pi ||y||^2_D}}$$
\end{Proposition}

{\bf Proof.} This is just the definition of the quotient from section 3, Definition 3.2.  \hfill \checkbox

\begin{Theorem} Suppose $\Delta$ is the absolute value of the discriminant of the number field $F.$ Then the first cohomology of an Arakelov divisor $D$ is the following ghost-space of the second kind.
$$H^1(D)=((I\otimes \R) /I)^{\omega},$$
where
 $$\omega = \frac{\sqrt{\Delta}}{e^{\deg D}}\cdot \sum \limits_{y \in I} e^{-\pi ||x+y||^2_D} \cdot m,$$
where $m$ is the Haar probability measure on $(I\otimes \R) /I$.
\end{Theorem}

{\bf Proof.} Obviously $\omega$ should be proportional to $\sum \limits_{y \in I} e^{-\pi ||x+y||^2_D} \cdot m.$ We just have to scale it to make it a probability measure. We have the following.
$$\int \limits_{\overline{x} \in (I\otimes \R) /I} \sum \limits_{y \in I} e^{-\pi ||x+y||^2_D} \cdot dm(\overline{x})= \int \limits_{x \in I \otimes \R} e^{-\pi ||x||^2_D} d M(x),$$
where $M$ is the measure on $I \otimes \R$ such that $I$ has covolume $1.$ 
If $M_D$ is the measure that corresponds to the hermitian metric $D,$ the above integral is equal to 
$$\frac{e^{\deg D}}{\sqrt{\Delta}}\cdot \int \limits_{x \in I \otimes \R} e^{-\pi ||x||^2_D} d M_D(x) $$
Now we just need to show that
$$\int \limits_{x \in I \otimes \R} e^{-\pi ||x||^2_D} d M_D(x) =1.$$
This is a pretty standard calculation. It can be done, e.g. by splitting up into the pieces that correspond to the infinite places of $F$ and using the following two identities.

1) (real factor) $$\alpha \int \limits_{x \in \R} e^{-\pi \alpha^2 x^2} dx =1$$

2)  (complex factor) $$\alpha \!\!\! \int \limits_{x+iy \in \C} e^{-\pi \alpha (x^2+y^2)} dx dy =1$$
These are very standard identities. The second one follows from the direct calculation in polar coordinates. The first one is essentially the square root of the second one.  \hfill \checkbox

\vskip 0.2cm

Now we are ready for the Serre's duality theorem. For this we  need to recall the definition of the canonical Arakelov divisor $K$ on $F.$ It is defined (cf., e.g.  \cite{GS}) as having associated fractional ideal $\partial ^{-1}$ and zero infinite components. Here $\partial $ is the different of $F.$

\begin{Theorem} (Serre's duality) For any Arakelov divisor $D$ we have the following duality of ghost-spaces.
$$H^1(D)=\widehat{H^0(K-D)}$$
\end{Theorem}

{\bf Proof.} First we need  to establish duality on the level of underlining locally compact groups. Suppose $I$ is the fractional ideal associated with $D.$ It follows from the definition of $K$ that $(I\otimes \R) /I \! = \! (F\otimes \R) /I$ is dual to $\partial ^{-1} I^{-1},$ where $\partial$ is the different of $F.$ The duality is given by the following pairing ($\overline{x} \in (F\otimes \R) /I, y\in \partial ^{-1} I^{-1}$).
$$(\overline{x},y)=e^{2\pi i \Tr (xy)},$$
where $x \in F\otimes \R $ is some representative of $\overline{x}$ and $\Tr (xy)$ is taken in the algebra $F\otimes \R.$

Now in order to prove the theorem we just need to show that for every $y\in \partial ^{-1} I^{-1}$
$$e^{-\pi ||y||^2_{K-D}} = \int \limits_{\overline{x} \in (I\otimes \R) /I} e^{2\pi i \Tr (x y)} d\omega (\overline{x}),$$
where $\omega $ is the probability measure from Theorem 5.1. Let's just simplify the right hand side.
$$ \int \limits_{\overline{x} \in (I\otimes \R) /I} \!\!\!\!\! \frac{\sqrt{\Delta}}{e^{\deg D}} \sum \limits_{z \in I} e^{-\pi ||x+z||^2_D} e^{2\pi i \Tr (xy)} dm(\overline{x}) \! = \!\!\!\!\!
\int \limits_{x \in I\otimes \R} \!\!\!\! e^{-\pi ||x||^2_D} e^{2\pi i \Tr (xy)} dM_D(x)$$  
This is a pretty standard integral. For the convenience of a reader, we reproduce the calculations in some details below.

Let us suppose that the infinite part of $D$ is given by the real numbers $(\sigma_1,...\sigma_{r_1},\sigma_{r_1+1}, ...\sigma_{r_1+r_2}).$ Splitting up the above integral, and $e^{-\pi ||y||^2_{K-D}}$ into the product of $r_1+r_2$ factors corresponding to different $\sigma_i,$ it is enough to prove the following two lemmas.

\begin{Lemma} (real factor) For any real $\sigma$ and $y$ the following identity is true.
$$\int\limits_{x\in \R} e^{-\pi e^{-2\sigma} x^2+2\pi ixy}\cdot e^{-\sigma} dx=e^{-\pi e^{2\sigma} y^2}$$
\end{Lemma}
{\bf Proof.} First of all, multiplying $x$ by $e^{-\sigma}$ and $y$ by $e^{\sigma }$ we can get rid of $\sigma .$ So we just need to prove that
$$\int\limits_{x\in \R} e^{-\pi  x^2+2\pi ixy}\cdot dx=e^{-\pi  y^2}.$$
The left hand side can be rewritten as 
$$\int\limits_{x\in \R} e^{-\pi (x+iy)^2}\cdot e^{-\pi y^2} dx$$
By contour integration, it is equal to
$$\int\limits_{x\in \R} e^{-\pi x^2}\cdot e^{-\pi y^2} dx = e^{-\pi y^2},$$
the lemma is proven.  \hfill \checkbox

\begin{Lemma} (complex factor) For any $\sigma \in R$ and $y=y_1+iy_2 \in \C$ the following identity is true.
$$ \int\limits_{x_1+ix_2 \in \C} e^{-2\pi e^{-\sigma} (x_1^2+x_2^2)}e^{4\pi i(x_1y_1-x_2y_2)}\cdot 2e^{-\sigma } dx_1 dx_2 = e^{-\pi \cdot 2e^{\sigma}(y_1^2+y_2^2)}$$
\end{Lemma}

{\bf Proof.} First of all, multiplying $x_1$ and $x_2$ by $e^{-\sigma /2},$ and $y_1$ and $y_2$ by $e^{-\sigma /2},$ we can get rid of $\sigma .$ 
So we just need to prove that
$$\int\limits_{x_1+ix_2 \in \C} e^{-2\pi  (x_1^2+x_2^2)} e^{4\pi i(x_1y_1-x_2y_2)}\cdot 2 dx_1 dx_2 = e^{-\pi \cdot 2(y_1^2+y_2^2)}$$
The left hand side can be rewritten as 
$$\int\limits_{x_1} \int\limits_{x_2} 2e^{-2\pi (x_1-iy_1)^2 -2\pi (x_2+iy_2)^2}\cdot e^{-2\pi (y_1^2+y_2^2)} dx_1 dx_2$$
This is equal to $e^{-2\pi (y_1^2+y_2^2)}$ by splitting up the above integral and then proceeding like in the previous lemma.  \hfill \checkbox

\vskip 0.2cm

So, we established the Serre's duality as the duality of ghost-spaces. The obvious corollary of it, and Theorem 4.1 is the following. 

\begin{Corollary} In the above notations,
$$h^1(D)=h^0(K-D)$$
\end{Corollary}

Now we obtain the Riemann-Roch formula using the additivity of dimension in the short exact sequences of ghost-spaces from section 3.

\begin{Theorem} (Riemann-Roch formula)
$$h^0(D)-h^1(D)=\deg D -\frac12 \log \Delta$$
\end{Theorem}

{\bf Proof.} We use the notations of Theorem 5.1. By Proposition 3.1 and Remark 3.2, 
$$h^1(D)=\dim H^1(D)=\dim_m (I\otimes \R)/I -\dim_m ((I\otimes \R)/I)_v=$$
$$= -\dim_m ((I\otimes \R)/I)_v  =   -(\dim_M(I\otimes \R)_u -\dim I_u)= $$
$$=  h^0(D)-\dim_M (I\otimes \R)_u $$
So we have that
$$h^0(D)-h^1(D)=\dim_M (I\otimes \R)_u =\log \int\limits_{x\in I\otimes R}e^{-\pi ||x||^2_D} dM_(x)=\log \frac{e^{\deg D}}{\sqrt{\Delta }}$$
as in the proof of Theorem 5.1. This proves the theorem.  \hfill \checkbox

So, we recovered the Riemann-Roch theorem of van der Geer and Schoof (first proven by Tate in his thesis). Our approach, of course, gives much more structure. We should also note that instead of using the Poisson summation formula, we basically reproved it along the lines of the usual proof of the Riemann-Roch theorem in the geometric case.

%%%%%%%%%%%%%%%%%%%%%%%%%%%%%%%%%%%%%%%%%%%
\section{Further remarks and open problems}

There are many directions in which the theory can be developed further. We list below the most interesting possibilities.

1) We believe that the theory can be extended to the higher-dimen\-sional case, at least to the case of curves over number fields. There we have $H^0(D),$ $H^1(D),$ and $H^2(D).$ We believe that $H^0(D)$ should be a discrete finite-dimensional ghost-space of the first kind. $H^2(D)$ should be a compact ghost-space of the second kind, dual to $H^0(K-D).$ The most troublesome part is $H^1(D).$ If $D$ has geometric degree at least $2g-1$ (for the curves of genus $g$) then $H^2(D)$ should be trivial, and $H^1(D)$ should be a compact ghost-space of the second kind. If $D$ has negative geometric degree then $H^0(D)$ is trivial, and $H^1(D)$ is a discrete ghost-space of the first kind. However the most interesting case of $0\le \deg D \le 2g-2$ is not covered above. In this case we conjecture that there still exists a ghost-space interpretation of $H^1(D),$ which is a locally compact group with the convolution structure that generalizes the structures of the ghost-spaces of the first and second kind as in is the following example.

{\bf Example.} Suppose $G$ is a locally compact abelian group, $u$ is an even continuous function on it, such that $u(0)=1$. Suppose also that $\mu $ is an even probability measure on $G$. Then the following convolution structure is commutative and associative.
$$\delta_x *\delta_y = \frac{u(x)u(y)}{u(x+y)}T_{x+y}\mu$$

This higher-dimensional generalization is clearly very important. Ultimately, one would like to translate from geometry such things as Kodaira-Spencer map to get a shot at the $abc$-type results. 

As a first step towards this goal, one should try to develop a ghost-space cohomology theory of hermitian line bundles on complex curves. This is still a one-dimensional problem but now all valuations are Archimedean. Possibly the underlining abelian groups here will be some functions spaces and will no longer be locally compact.
2) It is of some interest to extend the theory from the Arakelov divisors to the more general ``coherent ghost-sheaves", whatever this should mean. As a first step in this direction, Ichiro Miyada extended our theory to the higher rank locally free sheaves (cf. \cite{Miyada}). He also proposed a more adelic version of the theory.

3) As noted in \cite{GS}, prop. 6, zeta function of $F$ is kind of given by the following integral.
$$\int\limits_{\rm{Pic}(F)} e^{sh^0(D)+(1-s)h^1(D)}d[D]$$
In particular, Riemann zeta function is related to the family of ghost-spaces $\Z_u,$ where $u(x)=e^{-\pi \alpha x^2}$ for positive $\alpha .$ This extra structure of the ghost-space could be of some interest, as it relates arithmetic to harmonic analysis, which is coherent with some of the recent approaches to the Riemann Hypothesis.  We leave it to the RH specialists to figure out if it could be of any use.

4) The abstract theory of ghost-spaces, especially its analytic aspects are yet to be fully developed. First of all, one would like to develop  the theory of ``mixed ghost-spaces" i.e. groups with the convolution structures like in the Example above. One would also like to have a theory which is more symmetric with respect to duality. There are some obstacles here, some of which were resolved in \cite{ppd}.


\begin{thebibliography}{99}

\bibitem{BergForst} Christian Berg, Gunnar Forst. Potential theory on locally compact abelian groups. {\it Ergebnisse der Mathematik und ihrer Grenzgebiete, Band 87.} Springer-Verlag, New York-Heidelberg, 1975.

\bibitem{ppd} Alexandr Borisov. Positive positive-definite functions and measures on locally compact abelian groups, preprint (1999). Web address: http://www.math.psu.edu/borisov/

\bibitem{Folland} Gerald B. Folland. A course in abstract harmonic analysis. {\it Studies in Advanced Mathematics} CRC Press, 1995.
 
\bibitem{GS} Gerard van der Geer; Ren\'e Schoof. Effectivity of Arakelov divisors and the theta divisor of a number field, {\it Selecta Math.,} to appear. The preprint is available at
http://xxx.lanl.gov/abs/math/9802121.

\bibitem{Miyada} Ichiro Miyada.  Borisov-Serre duality of adelic spaces with convolution structure of measures
and Riemann-Roch formula for matrix divisors over a global field, preprint (1998).

\bibitem{Pym} John S. Pym. Weakly separately continuous measure algebras. {\it Math. Ann.} {\bf 175} (1968), 207--219.

\bibitem{Roesler1} Margit R\"osler. Convolution algebras which are not necessarily positivity-preserving. {\it Applications of hypergroups and related measure algebras (Seattle, WA 1993)}, 299--318, {\it Contemp. Math.} {\bf 183} AMS, Providence, RI, 1995.

\bibitem{Roesler2} Margit R\"osler. On the dual of commutative signed hypergroup. {\it Manuscripta Math.} {\bf 88} (1995), no. 2, 147--163.

\bibitem{Szpiro} Lucien Szpiro. Presentation de la theorie d'Arakelov.  (English) [Presentation of Arakelov theory] {\it Current trends in arithmetical algebraic geometry (Arcata, CA, 1985),} 279--293, {Contemp. Math.} {\bf 67} AMS, Providence, RI, 1987. 

\bibitem{Tate} J. T. Tate. Fourier analysis in number fields, and Hecke's zeta functions. 1967 {\it Algebraic Number Theory (Proc. Instructional Conf. Brighton, 1965)} pp. 305-347 {\it Thompson, Washington, D.C.}

\bibitem{Vain} Leonid I. Vainerman. The duality of algebras with involution and generalized shift operators. (Russian) Translated in {\it J. Soviet Math.} {\bf 42} (1988), no. 6, 2113--2137. {\it Itogi nauki i techniki, Mathematical analysis, } Vol. {\bf 24} (Russian), 165--205. Akad. Nauk Ssr, Vsesoyuz. Inst. Nauchn. i Techn. Inform., Moscow, 1986.

\bibitem{Voit} Michael Voit. A positivity result and normalization of positive convolution structures. {\it Math. Ann.} {\bf 297} (1993), no. 4, 677--692.



\end{thebibliography}
\end{document}